\documentclass{amsart}

\usepackage{amsthm}
\usepackage{amsmath}
\usepackage{amssymb}
\usepackage{amsfonts}
\usepackage{latexsym}
\usepackage{enumitem}
\usepackage{mathrsfs}
\usepackage[all]{xy} \SelectTips{eu}{}
\usepackage{hyperref}

\newtheorem{intthm}{Theorem}[]

\newtheorem*{intque*}{Question}

\newcommand{\numberseries}{\bfseries}   
\newlength{\thmtopspace}                
\newlength{\thmbotspace}                
\newlength{\thmheadspace}               
\newlength{\thmindent}                  
\newtheoremstyle{bfupright head,slanted body}
                {\thmtopspace}{\thmbotspace}
                {\slshape}{\thmindent}{\bfseries}{.}{\thmheadspace}
                {{\numberseries \thmnumber{#2\;}}\thmnote{#3}}

\newtheoremstyle{bfupright head,upright body}
                {\thmtopspace}{\thmbotspace}
                {\upshape}{\thmindent}{\bfseries}{.}{\thmheadspace}
                {{\numberseries \thmnumber{#2\;}}\thmnote{#3}}

\newtheoremstyle{fixed bf head,slanted body}
                {\thmtopspace}{\thmbotspace}{\slshape}
                {\thmindent}{\bfseries}{.}{\thmheadspace}
                {{\numberseries \thmnumber{#2\;}}\thmname{#1}\thmnote{ (#3)}}

\newtheoremstyle{fixed bf head,upright body}
                {\thmtopspace}{\thmbotspace}{\upshape}
                {\thmindent}{\bfseries}{.}{\thmheadspace}
                {{\numberseries \thmnumber{#2\;}}\thmname{#1}\thmnote{ (#3)}}

\newtheoremstyle{numbered paragraph}
                {\thmtopspace}{\thmbotspace}{\upshape}
                {\thmindent}{\upshape}{}{\thmheadspace}
                {{\numberseries \thmnumber{#2.}}}

\theoremstyle{bfupright head,slanted body}
\newtheorem{res}{}[section]             \newtheorem*{res*}{}

\theoremstyle{bfupright head,upright body}
\newtheorem{bfhpg}[res]{}               \newtheorem*{bfhpg*}{}

\theoremstyle{fixed bf head,slanted body}
          \newtheorem*{thm*}{Theorem}
\newtheorem{prp}[res]{Proposition}      \newtheorem*{prp*}{Proposition}
\newtheorem{cor}[res]{Corollary}        \newtheorem*{cor*}{Corollary}
\newtheorem{lem}[res]{Lemma}            \newtheorem*{lem*}{Lemma}
         \newtheorem*{que*}{Question}
\theoremstyle{fixed bf head,upright body}
       \newtheorem*{dfn*}{Definition}
\newtheorem{rmk}[res]{Remark}           \newtheorem*{rmk*}{Remark}
           \newtheorem*{exa*}{Example}

\theoremstyle{numbered paragraph}
\newtheorem{ipg}[res]{}

\newlength{\thmlistleft}        
\newlength{\thmlistright}       
\newlength{\thmlistpartopsep}   
\newlength{\thmlisttopsep}      
\newlength{\thmlistparsep}      
\newlength{\thmlistitemsep}     

\newcounter{eqc}
  {\end{list}}%

\newcounter{prt}
\newenvironment{prt}{\begin{list}{\upshape (\alph{prt})}%
    {\usecounter{prt}%
      \setlength{\leftmargin}{\thmlistleft}%
      \setlength{\labelwidth}{\thmlistleft}%
      \setlength{\rightmargin}{\thmlistright}%
      \setlength{\partopsep}{\thmlistpartopsep}%
      \setlength{\topsep}{\thmlisttopsep}%
      \setlength{\parsep}{\thmlistparsep}%
      \setlength{\itemsep}{\thmlistitemsep}}}%
  {\end{list}}%

\newcounter{rqm}
\newenvironment{rqm}{\begin{list}{\upshape (\arabic{rqm})}%
    {\usecounter{rqm}%
      \setlength{\leftmargin}{\thmlistleft}%
      \setlength{\labelwidth}{\thmlistleft}%
      \setlength{\rightmargin}{\thmlistright}%
      \setlength{\partopsep}{\thmlistpartopsep}%
      \setlength{\topsep}{\thmlisttopsep}%
      \setlength{\parsep}{\thmlistparsep}%
      \setlength{\itemsep}{\thmlistitemsep}}}%
  {\end{list}}%

\newenvironment{itemlist}{\nopagebreak \begin{list}{$\bullet$}%
    {\setlength{\leftmargin}{1.5em}%
      \setlength{\labelwidth}{\thmlistleft}%
      \setlength{\rightmargin}{\thmlistright}%
      \setlength{\partopsep}{\thmlistpartopsep}%
      \setlength{\topsep}{\thmlisttopsep}%
      \setlength{\parsep}{\thmlistparsep}%
      \setlength{\itemsep}{\thmlistitemsep}}}%
  {\end{list}}%

\newenvironment{prf*}[1][Proof]{%
  \begin{proof}[\bf #1]
    \setcounter{equation}{0}
    }
  {\end{proof}
}

\newcommand{\pgref}[1]{\ref{#1}}

\renewcommand{\eqref}[1]{(\pgref{eq:#1})}

\newcommand{\thmcite}[2][?]{\cite[thm.~#1]{#2}}
\newcommand{\prpcite}[2][?]{\cite[prop.~#1]{#2}}
\newcommand{\quecite}[2][?]{\cite[ques.~#1]{#2}}
\newcommand{\corcite}[2][?]{\cite[cor.~#1]{#2}}
\newcommand{\lemcite}[2][?]{\cite[lem.~#1]{#2}}

\newcommand{\dfncite}[2][?]{\cite[def.~#1]{#2}}
\newcommand{\rmkcite}[2][?]{\cite[rmk.~#1]{#2}}

\numberwithin{equation}{res}

\def\urltilda{\kern -.15em\lower .7ex\hbox{\~{}}\kern .04em}

\newcommand{\A}{\mathsf{A}}
\newcommand{\Smd}{\mathsf{Smd}}

\newcommand{\Mod}[1]{\mathsf{Mod}(#1)}

\newcommand{\xra}[2][]{\xrightarrow[#1]{\:#2\:}}

\newcommand{\Qop}{Q^{\sf op}}

\newcommand{\C}{\mbox{\rm C}}
\newcommand{\coker}{\mbox{\rm Coker}}

\newcommand{\Coker}[1]{\nobreak{\operatorname{Coker}#1}}
\newcommand{\Rep}[2]{{\sf Rep}(#1,#2)}

\def\bE{{\boldsymbol E}}

   \def\soft#1{\leavevmode\setbox0=\hbox{h}\dimen7=\ht0\advance
    \dimen7 by-1ex\relax\if t#1\relax\rlap{\raise.6\dimen7
    \hbox{\kern.3ex\char'47}}#1\relax\else\if T#1\relax
    \rlap{\raise.5\dimen7\hbox{\kern1.3ex\char'47}}#1\relax
    \else\if d#1\relax\rlap{\raise.5\dimen7\hbox{\kern.9ex
    \char'47}}#1\relax\else\if D#1\relax\rlap{\raise.5\dimen7
    \hbox{\kern1.4ex\char'47}}#1\relax\else\if l#1\relax
    \rlap{\raise.5\dimen7\hbox{\kern.4ex\char'47}}#1\relax
    \else\if L#1\relax\rlap{\raise.5\dimen7\hbox{\kern.7ex
    \char'47}}#1\relax\else\message{accent \string\soft
    \space #1 not defined!}#1\relax\fi\fi\fi\fi\fi\fi}

\begin{document}

\title[Completeness of induced cotorsion pairs in representation categories]%
{Completeness of induced cotorsion pairs in representation categories of rooted quivers}

\author[Z.X. Di]{Zhenxing Di}

\address{Zhenxing Di: Department of Mathematics, Huaqiao University, Quanzhou 362021, China}

\email{dizhenxing@163.com}

\author[L.P. Li]{Liping Li}

\address{Liping Li: Department of Mathematics, Hunan Normal University, Changsha 410081, China}

\email{lipingli@hunnu.edu.cn}

\author[L. Liang]{Li Liang}

\address{Li Liang: 1. Department of Mathematics, Shantou University, Shantou 515063, China
\ \ \ \ \ \ 2. School of Mathematics and Physics, Lanzhou Jiaotong University, Lanzhou 730070, China}

\email{lliangnju@gmail.com}

\urladdr{https://sites.google.com/site/lliangnju}

\author[F. Xu]{Fei Xu}

\address{Fei Xu: Department of Mathematics, Shantou University, Shantou 515063, China}

\email{fxu@stu.edu.cn}

\thanks{Z.X. Di was partly supported by NSF of China grant 11971388;
L.P. Li was partly supported by NSF of China grant 11771135,
the Hunan Provincial Science and Technology Department grant 2019RS1039,
and Hunan Provincial Education Department grant 18A016;
L. Liang was partly supported by NSF of China grant 11761045;
F. Xu was partly supported by NSF of China grant 11671245.}

\date{\today}



\begin{abstract}
This paper focuses on a question raised by Holm and J{\o}rgensen, who asked if the induced cotorsion pairs $(\Phi({\sf X}),\Phi({\sf X})^{\perp})$ and $(^{\perp}\Psi({\sf Y}),\Psi({\sf Y}))$ in $\Rep{Q}{\sf{A}}$--the category of all $\A$-valued representations of a quiver $Q$--are complete whenever $(\sf X,\sf Y)$ is a complete cotorsion pair in an abelian category $\sf{A}$ satisfying some mild conditions. Recently, Odaba\c{s}{\i} gave an affirmative answer if the quiver $Q$ is rooted and the cotorsion pair $(\sf X,\sf Y)$ is further \emph{hereditary}. In this paper, we improve Odaba\c{s}{\i}'s work by removing the hereditary assumption on the cotorsion pair. As an application,
we show under certain mild conditions that if a subcategory $\sf L$,
which is not necessarily closed under direct summands,
of $\sf A$ is special precovering (resp., preenveloping),
then $\Phi(\sf L)$ (resp., $\Psi(\sf L)$) is special precovering (resp., preenveloping) in $\Rep{Q}{\A}$.
\end{abstract}

\maketitle

\thispagestyle{empty}
\section*{Introduction}
\noindent
Holm and J{\o}rgensen \cite{HJ19} proved that if an abelian category $\sf{A}$ satisfies some mild conditions, then every cotorsion pair $(\sf X,\sf Y)$ in $\sf{A}$ can induce the cotorsion pairs
\begin{center}
$(\Phi({\sf X}),\Phi({\sf X})^{\perp})$ \quad and \quad $(^{\perp}\Psi({\sf Y}),\Psi({\sf Y}))$
\end{center}
in $\Rep{Q}{\sf{A}}$, the category of all $\A$-valued representations of a quiver $Q$. Here $\Phi({\sf X})$ and $\Psi({\sf Y})$ are defined as follows (see \ref{quiver notation} for details):
\begin{itemlist}
\item $\Phi({\sf X})=\left\{ X\in\Rep{Q}{\A} \:
      \left|
        \begin{array}{c}
          \varphi_i^X\ \text {is a monomorphism and}\\
          \C_i(X)\in{\sf X}\ \text{for each vertex}\ i\ \text{in}\ Q
        \end{array}
      \right.
    \right\}$, and
\item $\Psi({\sf Y})=\left\{ Y\in\Rep{Q}{\A} \:
      \left|
        \begin{array}{c}
          \psi_i^Y\ \text {is an epimorphism and}\\
          \mathrm{K}_i(Y)\in{\sf Y}\ \text{for each vertex}\ i\ \text{in}\ Q
        \end{array}
      \right.
    \right\}$.
\end{itemlist}
This result gives a ``quiver-representation version" of \prpcite[3.6]{Gi04}
by Gillespie on cotorsion pairs in the category of chain complexes. At the same time, Holm and J{\o}rgensen raised in \quecite[7.7]{HJ19} the following question:
\begin{intque*}[Holm and J{\o}rgensen]
Is it true that if the cotorsion pair $(\sf X,\sf Y)$ in $\sf{A}$ is complete then so are the above induced cotorsion pairs\,?
\end{intque*}

It seems hard to give an affirmative answer for any quiver $Q$. However, inspired by the work of Yang and Ding \cite{YD15} (see also Yang and Liu \cite{Yang-Liu}), Odaba\c{s}{\i} proved under some mild conditions that if the quiver $Q$ is left (resp., right) rooted and a complete cotorsion pair $(\sf X,\sf Y)$ in $\sf{A}$ is further \emph{hereditary}, then the induced cotorsion pair $(\Phi(\sf X),\Phi(\sf X)^\perp)$ (resp., $({^\perp\Psi(\sf Y)},\Psi(\sf Y))$) in $\Rep{Q}{\A}$ is complete as well; see \thmcite[5.6]{Od19}. In the present paper, we improve the above Odaba\c{s}{\i}'s results by removing the condition that the given cotorsion pair in $\A$ is hereditary.

\begin{intthm}\label{AAA}
Let $\sf A$ be an abelian category, and let $(\sf X,\sf Y)$ be a complete cotorsion pair in $\sf A$.
\begin{prt}
\item Suppose that $Q$ is a left rooted quiver and $\sf A$ has exact $|Q_1^{\ast\to i}|$-indexed coproducts for every vertex $i$ in $Q$ (e.g., $\A$ satisfies the axiom AB4). Then $(\Phi(\sf X),\Phi(\sf X)^\perp)$ is a complete cotorsion pair in $\Rep{Q}{\A}$.
\item Suppose that $Q$ is a right rooted quiver and $\sf A$ has exact $|Q_1^{i\to \ast}|$-indexed products for every vertex $i$ in $Q$ (e.g., $\A$ satisfies the axiom AB4*). Then $({^\perp\Psi(\sf Y)},\Psi(\sf Y))$ is a complete cotorsion pair in $\Rep{Q}{\A}$.
\end{prt}
\end{intthm}

Notice that there do exist a large number of cotorsion pairs
which is complete but not hereditary.
For example, it is known from \thmcite[4.1.6]{RT12} that over an arbitrary associated ring $R$,
the pair $({^\perp\sf Fp},{\sf Fp})$ forms a complete cotorsion pair in the category $\Mod{R}$ of left $R$-modules, where
${\sf Fp}$ denotes the subcategory of $\Mod{R}$
consisting of all FP-injective left $R$-modules.
However, the cotorsion pair $({^\perp\sf Fp},{\sf Fp})$ is hereditary if and only if $R$ is left coherent; see \thmcite[5.5]{BD-PM}.
\begin{center}
  $\ast \ \ \ast \ \ \ast$
\end{center}

\noindent
Originated from the concept of injective envelopes, the approximation theory has attracted increasing interest, and hence obtained considerable development especially in the context of abelian categories; see for example Auslander and Reiten \cite{MAsIRt91}, and Enochs and Jenda \cite{rha}. The notions of special precovering and preenveloping subcategories are important cornerstones in the approximation theory. With these notions in hand, one can obtain ``proper" resolutions under which the corresponding derived functors can be defined.

An important approach for producing special precovering subcategories is via complete cotorsion pairs. According to Hu and Zhu \lemcite[3.9]{Hu-Zhu},
one gets that a subcategory $\sf L$ in an abelian category with enough injectives is special precovering and closed under direct summands if and only if
$(\sf L,\sf L^\perp)$ forms a complete cotorsion pair. However, there do exist special precovering subcategories that may not be closed under direct summands, and example of which is the subcategory of $\Mod{R}$ consisting of all free $R$-modules.

As an application of Theorem \ref{AAA},
we obtain the following result.
It asserts under some mild conditions that
if the quiver $Q$ is left (resp., right) rooted,
then any special precovering (resp., preenveloping) subcategory,
which is not necessarily closed under direct summands, in $\sf A$
can induce a special precovering (resp., preenveloping) subcategory in $\Rep{Q}{\A}$.

\begin{intthm}\label{BBBB}
Let $\sf A$ be an abelian category, and let $\sf L$ be a subcategory of $\sf A$.
\begin{prt}
\item Suppose that $Q$ is a left rooted quiver, and
      suppose that $\sf A$ satisfies the axiom AB4 and
      has enough injectives.
      If $\sf L$ is special precovering in $\sf A$,
      then $\Phi(\sf L)$ is special precovering in $\Rep{Q}{\A}$.
\item Suppose that $Q$ is a right rooted quiver, and
      suppose that $\sf A$ satisfies the axiom AB4* and
      has enough projectives.
      If $\sf L$ is special preenveloping in $\sf A$,
      then $\Psi(\sf L)$ is special preenveloping in $\Rep{Q}{\A}$.
\end{prt}
\end{intthm}

The paper is organized as follows.
Section \ref{pre} contains some necessary notation and terminologies for use throughout this paper.
In Section \ref{tp}, we mainly give the proof of Theorem \ref{AAA}.
In Section \ref{model}, we give the proof of Theorem \ref{BBBB}.
Finally, we give an Appendix for reproving Theorem \ref{AAA}(a)
over a more simple quiver with 4 vertices
to comprehend the idea of the proof.

\section{Preliminaries}\label{pre}
\noindent
In this section we mainly recall some necessary notions and definitions.
\textsf{Throughout the paper},
let $\sf A$ be an abelian category and $Q=(Q_0,Q_1)$ a quiver.
In what follows,
by the term \emph{``subcategory''} we always mean a full additive subcategory closed under isomorphisms.

\begin{ipg}\label{AB}
Recall that an abelian category satisfies the axiom $AB3$
if it has small coproducts, equivalently, if it is cocomplete.
It satisfies the axiom $AB4$ if it satisfies the axiom $AB3$ and any coproduct of monomorphisms is a monomorphism.
The axioms $AB3^*$ and $AB4^*$ are dual to the axioms $AB3$ and $AB4$, respectively.
\end{ipg}

\begin{bfhpg}[\bf Quivers and representations]\label{quiver notation1}
A quiver $Q=(Q_0,Q_1)$ is actually a directed graph with vertex set $Q_0$ and arrow set $Q_1$. The symbol $\Qop$ denotes the opposite quiver of $Q$, which is obtained by reversing arrows of $Q$. For an arrow $a\in Q_1$, we always write $s(a)$ for its source and $t(a)$ for its target. For a vertex $i\in Q_0$ we let $Q_1^{i\to\ast}$ denote the set $Q_1^{i\to\ast}=\{a\in Q_1\,|\,s(a)=i\}$, and let $Q_1^{\ast\to i}$ denote the set $Q_1^{\ast\to i}=\{a\in Q_1\,|\,t(a)=i\}$.

By a \emph{representation} $X$ of $Q$ we mean a functor from $Q$ to $\A$,
which is determined by associating an object in $\A$ to each vertex $i\in Q_0$ and
a morphism $X(a): X(i)\to X(j)$ to each arrow $a: i\to j$ in $Q_1$. For representations $X$ and $Y$ of $Q$, a \emph{morphism} $f$ from $X$ to $Y$ is a natural transformation,
that is, a family of morphisms $\{f(i):X(i)\to Y(i)\}_{i\in Q_0}$
such that $Y(a)\circ f(i)=f(j)\circ X(a)$ for each arrow $a: i\to j$. All representations of $Q$ constitute an abelian category, which is denoted $\Rep{Q}{\A}$.
\end{bfhpg}

\begin{ipg}
For every vertex $i \in Q_0$, there exists an \emph{evaluation functor}
$${\sf e}_i: \Rep{Q}{\A}\to \A,$$
which maps an $\A$-valued representation $X$ of $Q$ to the object ${\sf e}_i(X) = X(i) \in\A$.
On the other hand, for every $i \in Q_0$ there exists a \emph{stalk functor}
$${\sf s}_i: \A \to \Rep{Q}{\A},$$
which to an object $M \in \A$ assigns the stalk representation ${\sf s}_i(M)$
given by ${\sf s}_i(M)(i) = M$ and ${\sf s}_i(M)(j) = 0$ for $j \neq i$;
for every arrow $a \in Q_1$ the morphism ${\sf s}_i(M)(a)$ is zero.
\end{ipg}

\begin{ipg}
Let $S$ be a set.
Following the notation given in \cite[2.3]{Od19},
denote by $|S|$ the cardinality of the set $S$.
The category $\sf A$ is said to \emph{have $|Q_1^{\ast\to i}|$-indexed coproducts} for some vertex $i$ in $Q_0$
if any family $\{C_u\}_{u\in U}$ of objects in $\sf A$
has the coproduct in $\sf A$ whenever $|U| \leqslant |Q_1^{\ast\to i}|$.
Furthermore, if any $|Q_1^{\ast\to i}|$-indexed coproduct of monomorphisms is a monomorphism,
we say that $\sf A$ has \emph{exact $|Q_1^{\ast\to i}|$-indexed coproducts}. The \emph{(exact) $|Q_1^{i\to\ast}|$-indexed products} can be defined dually.
\end{ipg}

\begin{ipg}\label{f and g}
Suppose that $\A$ has the coproduct of each family $\{C_u\}_{u\in U}$ where $|U|\leqslant|Q(i, j)|$ for any vertexes $i$ and $j$ in $Q_0$.
For a given vertex $i\in Q_0$ and any object $M\in \A$,
according to \cite[3.1]{HJ19},
the representation ${\sf f}_i(M)\in \Rep{Q}{\A}$ is defined as follows:
For each $j\in Q_0$, set ${\sf f}_i(M)(j)=\coprod_{p\in Q(i,j)}M_p,$
where each $M_p$ is a copy of $M$;
if there are no paths in $Q$ from $i$ to $j$,
then this coproduct is empty and hence ${\sf f}_i(M)(j)$ is $0$ in $\A$.
For an arrow $a : j \to k \in Q_1$ and $p \in Q(i, j)$,
then $ap$ is a path from $i$ to $k$.
Therefore, define ${\sf f}_i(M)(a)$ to be the unique morphism in $\A$
that makes the following diagram commutative:
$$\xymatrix{M_p\ar@{=}[rr]\ar[d]& & M_{ap}\ar[d]\\
{\sf f}_i(M)(j)\ar[rr]^{{\sf f}_i(M)(a)}& &{\sf f}_i(M)(k)}$$
in which the vertical morphisms are the canonical injections.
It is evident that the assignment $M \to {\sf f}_i(M)$ yields a functor
${\sf f}_i : \A \to \Rep{Q}{\A}$.

Suppose that $\A$ has the product of each family $\{C_u\}_{u\in U}$
where $|U| \leqslant|Q(i, j)|$ for any vertexes $i$ and $j$ in $Q_0$.
Then the opposite category ${\sf A^{op}}$ satisfies the above conditions
for the opposite quiver $\Qop$.
Let ${\sf g}_i$ denote the functor $({\sf f}_i)^{\sf op}$.
\end{ipg}

\begin{ipg}\label{quiver notation}
For each representation $X\in\Rep{Q}{\A}$ and each vertex $i\in Q_0$,
by the universal property of coproducts,
there is a unique morphism $$\varphi_{i}^{X}:\oplus_{a\in Q_{1}^{\ast\to i}}X(s(a))\to X(i).$$
Let $\C_i(X)$ denote the cokernel of $\varphi_{i}^{X}$;
it yields a functor from $\Rep{Q}{\A}$ to $\A$.
Dually, there is a unique morphism
$$\psi_{i}^{X}:X(i)\to \Pi_{a\in Q_{1}^{i\to \ast}}X(t(a)).$$
The symbol $\mathrm{K}_i(X)$ denotes the kernel of $\psi_{i}^{X}$,
which also yields a functor from $\Rep{Q}{\A}$ to $\A$.

For a subcategory $\sf X$ of $\A$, we set

\begin{itemlist}
\item $\Rep{Q}{{\sf X}}=\{X\in\Rep{Q}{\A}~|~X(i)\in {\sf X} \ \text {for each vertex}\ i\in Q_0\}$,
\item $\Phi({\sf X})=\left\{ X\in\Rep{Q}{\A} \:
      \left|
        \begin{array}{c}
          \varphi_i^X\ \text {is a monomorphism and}\\
          \C_i(X)\in{\sf X}\ \text{for each vertex}\ i\in Q_0
        \end{array}
      \right.
    \right\}$ and
\item $\Psi({\sf X})=\left\{ X\in\Rep{Q}{\A} \:
      \left|
        \begin{array}{c}
          \psi_i^X\ \text {is an epimorphism and}\\
          \mathrm{K}_i(X)\in{\sf X}\ \text{for each vertex}\ i\in Q_0
        \end{array}
      \right.
    \right\}$
\end{itemlist}
\end{ipg}

\begin{lem}\label{cf=0}
Let $Q$ be a quiver without oriented cycles. Then for each vertex $i\in Q_0$ there are natural isomorphisms $\C_i{\sf f}_i\cong{\sf id}_\A$ and $\mathrm{K}_i{\sf g}_i\cong{\sf id}_\A$.
\end{lem}
\begin{prf*}
We only prove the first isomorphism in the statement.

Let $M$ be an object in $\A$. One has ${\sf f}_i(M)(s(a))=0$ for any arrow $a\in Q_1^{\ast\to i}$, as there exist no path from $i$ to $s(a)$.
This implies that $\varphi_i^{{\sf f}_i(M)}=0$ with $\C_i({\sf f}_i(M))\cong{\sf f}_i(M)(i)=M$. So one has $\C_i{\sf f}_i\cong{\sf id}_\A$; the naturality holds clearly.
\end{prf*}

\begin{bfhpg}[\bf Rooted quivers]\label{rooted}
According to \cite{EOT04},
there is a transfinite sequence $\{V_{\alpha}\}_{\alpha \mathrm{ordinal}}$ of subsets of $Q_0$ as follows:

For the first ordinal $\alpha=0$ set $V_0=\emptyset$, for a successor ordinal $\alpha+1$ set
$$V_{\alpha+1}=\{i\in Q_0~|~i\ \text{is not the target of any arrow}\ a\in Q_1\ \text{with}\ s(a)\notin\cup_{\beta\leq\alpha}V_{\beta}\},$$
and for a limit ordinal $\alpha$ set $V_{\alpha}=\cup_{\beta<\alpha}V_{\beta}$.

It is clear that $V_1=\{i\in Q_0~|~\text{there is no arrow}\ a\in Q_1\ \text{with}\ t(a)=i\}$.
By \lemcite[2.7]{HJ19} and \corcite[2.8]{HJ19}, there is a chain
$V_1\subseteq V_2\subseteq\cdots\subseteq Q_0,$ and
if $a: i\to j$ is an arrow in $Q_1$ with $j\in V_{\alpha+1}$ for some ordinal $\alpha$, then $i$ must be in $V_{\alpha}$.

Following \dfncite[3.5]{EOT04} a quiver $Q$ is called \emph{left rooted}
if there exists an ordinal $\lambda$ such that $V_{\lambda}=Q_0$. In a dual manner,
a quiver $Q$ is said to be \emph{right rooted} if the opposite quiver $\Qop$ is left rooted.
\end{bfhpg}

The next result is from \prpcite[3.6]{EOT04}, which asserts that left (right) rooted quivers constitute quite a large class of quivers.

\begin{lem}\label{e-rooted}
Let $Q$ be a quiver. Then the following statements hold.
\begin{prt}
\item $Q$ is left rooted if and only if
it has no infinite sequence of arrows of the form $\cdots \to \bullet\to \bullet \to \bullet$ (not necessarily different).
\item $Q$ is right rooted if and only if
it has no infinite sequence of arrows of the form $\bullet \to \bullet \to \bullet \to \cdots$ (not necessarily different).
\end{prt}
\end{lem}

The next result,
which is from \prpcite[7.2]{HJ19},
will be used frequently in the sequel.

\begin{lem}\label{X(i) in X}
Let $\sf X$ be a subcategory of $\sf A$ which is closed under extensions.
\begin{prt}
\item Suppose that $Q$ is a left rooted quiver and $\sf A$ has $|Q_1^{\ast\to i}|$-indexed coproducts for every vertex $i \in Q_0$. If $\sf X$ is closed under $|Q_1^{\ast\to i}|$-indexed coproducts for every vertex $i \in Q_0$,
      then every representation $X \in \Phi({\sf X})$ has values in $\sf X$;
      that is, $X(i) \in \sf X$ for each vertex $i \in Q_0$.
\item Suppose that $Q$ is a right rooted quiver and $\sf A$ has $|Q_1^{i\to\ast}|$-indexed products for every vertex $i \in Q_0$. If $\sf X$ is closed under $|Q_1^{i\to\ast}|$-indexed products for every vertex $i \in Q_0$,
      Then every representation $X \in \Psi({\sf X})$ has values in $\sf X$;
      that is, $X(i) \in \sf X$ for each vertex $i \in Q_0$.
\end{prt}
\end{lem}

\begin{ipg}\label{rooted}
Let $\sf X$ be a subcategory of $\A$.
Set
$${\sf X}^{\bot}=\{M\in {\sf A}\ |\ \textrm{Ext}_{\sf A}^{1}(X,M)=0 {\rm \ for\ each\ object}\ X\in{\sf X}\},\ \text{and}$$
$$^{\bot}{\sf X}=\{M\in {\sf A}\ |\ \textrm{Ext}_{\sf A}^{1}( M,X)=0 {\rm \ for\ each\ object}\ X\in{\sf X}\}.$$
A morphism $f : X \to M$ is called a \emph{special} $\sf X$-\emph{precover} of an object $M$
if $f$ is an epimorphism, $X\in \sf X$ and Ker$(f) \in \sf X^\perp$.
Dually, a morphism $g : M \to X$ is called a \emph{special} $\sf X$-\emph{preenvelope} of $M$
if $g$ is a monomorphism, $X \in\sf X$ and $\coker (g) \in {^\perp\sf X}$.

The subcategory $\sf X$ is called \emph{special precovering} (resp., \emph{special preenveloping})
if every object in $\A$ has a special $\sf X$-precover (resp., special $\sf X$-preenvelope).
\end{ipg}

\begin{lem}\label{an-side complete}
Let $\sf X$ be a subcategory of $\sf A$.
\begin{prt}
\item If $\sf X$ is special preenveloping,
      then each representation in $\Rep{Q}{\A}$ can be embedded in a representation in $\Rep{Q}{\sf X}$.
\item If $\sf X$ is special precovering,
      then each representation in $\Rep{Q}{\A}$ is an epimorphic image of a representation in $\Rep{Q}{\sf X}$.
\end{prt}
\end{lem}

\begin{prf*}
We only prove (a); one can prove (b) dually.

Let $M$ be a representation in $\Rep{Q}{\A}$.
Since $\sf X$ is special preenveloping in $\sf A$,
there exists for every vertex $i \in Q_0$
a short exact sequence
$$0\to M(i)\overset{f(i)} \longrightarrow X(i)\overset{g(i)} \longrightarrow Y(i)\to0$$
in $\sf A$ with $X(i)\in \sf X$ and $Y(i)\in \sf X^\perp$.
Hence, for any arrow $a:i\to j\in Q_1$
there exists a morphism $X(a)$ (hence, a morphism $Y(a)$) such that the following diagram
$$\xymatrix@C=22pt@R=22pt{
0\ar[r] &M(i)\ar[r]^{f(i)}\ar[d]_{M(a)}&X(i) \ar@.[d]_{X(a)} \ar[r]^{g(i)}& Y(i) \ar@.[d]_{Y(a)}\ar[r]&0  \\
0\ar[r] &M(j)\ar[r]^{f(j)}&  X(j)  \ar[r]^{g(j)} & Y(j)\ar[r]&0}$$
is commutative.
Thus, we obtain a short exact sequence
$$0\to M \to X\to Y\to0$$
in $\Rep{Q}{\A}$ such that $X\in \Rep{Q}{\sf X}$.
\end{prf*}

\begin{bfhpg}[\bf Cotorsion pairs]\label{cotpair}
Let $\sf X$ and $\sf Y$ be two subcategories of $\A$.
The pair $(\sf X,\sf Y)$ is called a \emph{cotorsion pair} \cite{rha}
if $\sf{X}^{\bot}=\sf{Y}$ and $^{\bot}\sf{Y}=\sf{X}$.

Let $(\sf X,\sf Y)$ be a cotorsion pair in $\sf A$.
Then $\sf A$ is said to \emph{have enough} $\sf X$-\emph{objects} (resp., $\sf Y$-\emph{objects})
if for every object $M$ in $\A$,
there exists an epimorphism $X\twoheadrightarrow M$ (resp., a monomorphism $M \rightarrowtail Y$) with $X \in \sf X$ (resp., $Y \in \sf Y$).
Following from \cite{rha} that a cotorsion pair $(\sf X,\sf Y)$ is said to \emph{have enough projectives}
(resp., \emph{injectives})
if every object of $\A$ has a special $\sf X$-precover (resp., a special $\sf Y$-preenvelope).
A cotorsion pair is called \emph{complete} if it has both enough projectives and enough injectives.
\end{bfhpg}

The following is a general form of the so-called Salce's Lemma (see \cite{S-L}),
which can be found in \lemcite[3.3]{Od19}.

\begin{lem}\label{complete}
Let $(\sf X, \sf Y)$ be a cotorsion pair in $\A$.
\begin{prt}
\item If $(\sf X, \sf Y)$ has enough projectives and $\A$ has enough $\sf Y$-objects
(e.g. $\A$ has enough injective objects),
then the cotorsion pair $(\sf X, \sf Y)$ is complete.
\item If $(\sf X, \sf Y)$ has enough injectives and $\A$ has enough $\sf X$-objects
(e.g. $\A$ has enough projective objects),
then the cotorsion pair $(\sf X, \sf Y)$ is complete.
\end{prt}
\end{lem}

\section{Completeness of the induced cotorsion pairs}\label{tp}
\noindent
In this section we mainly give the proof of Theorem \ref{AAA} advertised in the introduction. We begin with the following lemma, which was first proved by Holm and J{\o}rgensen \thmcite[7.4 and 7.9]{HJ19}; one refers to \rmkcite[3.12 and 3.13]{Od19} for the next version.

\begin{lem}\label{containment}
Let $(\sf X,\sf Y)$ be a cotorsion pair in $\sf A$ and $Q$ a quiver.
\begin{prt}
\item Suppose that $\sf A$ has enough $\sf{Y}$-objects, and has $|Q_1^{\ast\to i}|$-indexed coproducts for every vertex $i\in Q_0$. Then $(\Phi({\sf X}),\Phi({\sf X})^\perp)$ is a cotorsion pair in the category $\Rep{Q}{\sf A}$. If furthermore $Q$ is a left rooted quiver, then $\Rep{Q}{\sf Y}\subseteq \Phi({\sf X})^\perp$.
\item Suppose that $\sf A$ has enough $\sf{X}$-objects, and has $|Q_1^{i\to \ast}|$-indexed products for every vertex $i$ in $Q_0$. Then $(^\perp\Psi({\sf Y}),\Psi({\sf Y}))$ is a cotorsion pair in the category $\Rep{Q}{\sf A}$. If furthermore $Q$ is a right rooted quiver, then $\Rep{Q}{\sf X}\subseteq {^\perp\Psi({\sf Y})}$.
\end{prt}
\end{lem}

\begin{rmk}\label{equality}
Let $(\sf X,\sf Y)$ be a cotorsion pair in $\sf A$. If $\sf A$ satisfies the axiom AB4 and has enough $\sf{Y}$-objects, and $Q$ is a left rooted quiver, then by a similar proof as in \thmcite[7.9(a)]{HJ19} one gets that $\Rep{Q}{\sf Y}=\Phi({\sf X})^\perp$ holds true. Dually, if $\sf A$ satisfies the axiom AB4* and has enough $\sf{X}$-objects, and $Q$ is a right rooted quiver, then $\Rep{Q}{\sf X}={^\perp\Psi({\sf Y})}$.
\end{rmk}

The next result asserts that under some mild conditions
the converses of statements in Lemma \ref{containment} hold true.

\begin{prp}\label{another side of cp}
Let $\sf X$ be a subcategory of $\sf A$.
\begin{prt}
\item Suppose that $\sf A$ satisfies the axiom $AB4$ and $Q$ is a left rooted quiver. If $(\Phi({\sf X}),\Rep{Q}{\sf{X}^\perp})$ is a cotorsion pair in $\Rep{Q}{\A}$,
      then $(\sf X,\sf{X}^\perp)$ is a cotorsion pair in $\sf A$.
\item Suppose that $\sf A$ satisfies the axiom $AB4^*$ and $Q$ is a right rooted quiver. If $(\Rep{Q}{{^\perp\sf{X}}},\Psi({\sf X}))$ is a cotorsion pair in $\Rep{Q}{\A}$,
      then $({^\perp\sf{X}},\sf X)$ is a cotorsion pair in $\sf A$.
\end{prt}
\end{prp}
\begin{prf*}
We only give the proof of (a); the proof of (b) is dually.

It is sufficient to show that $^\perp(\sf X^\perp) \subseteq \sf X$.
Let $M$ be in $^\perp(\sf X^\perp)$, and $Y$ a representation in $\Rep{Q}{\sf{X}^\perp}$.
By \prpcite[5.2(a)]{HJ19}, for any vertex $i\in Q_0$,
one has
$$\textrm{Ext}_{\Rep{Q}{\A}}^1({\sf f}_i(M),Y)\cong \textrm{Ext}_{\A}^1(M,{\sf e}_i(Y))=0.$$
Hence, we conclude that ${\sf f}_i(M)\in \Phi({\sf X})$,
as $(\Phi({\sf X}),\Rep{Q}{\sf{X}^\perp})$ forms a cotorsion pair in $\Rep{Q}{\A}$ by assumption.
Note that $Q$ is left rooted.
It follows from Lemma \ref{e-rooted}(a) that $Q$ has no oriented cycles, and so $\C_i({\sf f}_i(M))\cong{\sf f}_i(M)(i)=M$ by Lemma \ref{cf=0}.
However, $\C_i({\sf f}_i(M))\in \sf X$, as we have proven that ${\sf f}_i(M)\in \Phi({\sf X})$.
It follows that $M\in\sf X$, as desired.
\end{prf*}

We now give the proof of Theorem \ref{AAA}; one refers to Appendix for a proof over a more simple quiver with 4 vertices to comprehend the idea of the proof.

\begin{bfhpg}[\bf Proof of Theorem \ref{AAA}]\label{PROOF OF a1}
We only prove (a); the statement (b) is proved dually.

Let $X$ be a representation in $\Rep{Q}{\A}$.
Since $(\sf X,\sf Y)$ is a complete cotorsion pair in $\sf A$,
one concludes that for each vertex $i\in Q_0$
there exists a short exact sequence
$$0\to B(i) \overset{k(i)}\longrightarrow A(i) \overset{h(i)}\longrightarrow X(i) \to 0$$
in $\A$ with $A(i)\in \sf X$ and $B(i)\in\sf Y$.
In particular, for each arrow $a: i\to j \in Q_1$ and for each morphism
$\xymatrix{A(i)\ar[rr]^{ X(a)\circ h(i)}& & X(j)}$
there exists a morphism $A(a):A(i)\to A(j)$ (and hence a morphism $B(a):B(i)\to B(j)$)
such that the diagram
$$\xymatrix{
0\ar[r] &B(i)\ar[r]^{k(i)}\ar@.[d]_{B(a)}& A(i) \ar@.[d]_{A(a)} \ar[r]^{h(i)}&  X(i) \ar[d]_{X(a)}\ar[r]&0  \\
0\ar[r] &B(j)\ar[r]^{k(j)}&  A(j)  \ar[r]^{h(j)} & X(j)\ar[r]&0}$$
is commutative. Hence one gets a short exact sequence
$$\bE:\,\,\,0\to B\to A\to X\to 0$$
in $\Rep{Q}{\A}$ such that $A\in \Rep{Q}{\sf X}$ and $B\in \Rep{Q}{\sf Y}$.
But $A$ may fail in $\Phi(\sf X)$.
We will use transfinite induction to construct a short exact sequence
$$\bE':\,\,\,0\to B'\to A'\to X\to 0$$
in $\Rep{Q}{\A}$ such that
$A'\in \Phi(\sf X)$ and $B'\in\Rep{Q}{\sf Y}\subseteq \Phi({\sf X})^\perp$ (see Lemma \ref{containment}(a)).
If we have done, then the cotorsion pair $(\Phi(\sf X),\Phi(\sf X)^\perp)$ is complete in $\Rep{Q}{\A}$ by Lemma \ref{complete}(a), as $\Rep{Q}{\A}$ has enough $\Phi(\sf X)^\perp$-objects by Lemmas \ref{an-side complete}(a) and \ref{containment}(a).

Let $\{V_{\alpha}\}$ be the transfinite sequence of subsets of $Q_0$. Since $Q$ is left rooted,
one has $Q_0=V_{\lambda}$ for some ordinal $\lambda$.
Next we construct a continuous inverse $\lambda$-sequence
$$\{\,\bE_\alpha:\,\,\, 0\to B_\alpha \overset{k_\alpha}\longrightarrow A_\alpha\overset{h_\alpha}\longrightarrow X_\alpha\to0 \,\}_{\alpha\leqslant\lambda}$$
of short exact sequences in $\Rep{Q}{\A}$ satisfying the following conditions:

\begin{prt}
\item For every ordinal $0<\alpha\leqslant\lambda$, $X_\alpha=X$.
\item For every ordinal $0<\alpha\leqslant\lambda$ and for each $i\in Q_0\backslash V_\alpha$,
$A_\alpha(i) = A(i)\in \sf{X}$ and $B_\alpha(i) = B(i)\in \sf{Y}$.
\item For every ordinal $\alpha\leqslant\lambda$ and each $i\in V_\alpha$,
$\varphi_i^{A_\alpha}$ is a monomorphism,
$A_\alpha(i)\in \sf{X}$, $\C_i(A_\alpha)\in \sf{X}$ and $B_\alpha(i)\in \sf{Y}$.
\item For every $0<\alpha<\alpha'\leqslant\lambda$,
$A_\alpha(i) = A_{\alpha'}(i)$ and $B_\alpha(i) = B_{\alpha'}(i)$
for all $i\notin V_{\alpha'}\backslash V_\alpha$,
and there exists the following commutative diagram in $\Rep{Q}{\A}$
$$\xymatrix{
\bE_{\alpha'}:\quad    0\ar[r] & B_{\alpha'}\ar[r]^{k_{\alpha'}}\ar[d]^{{g_{{\alpha},{\alpha'}}}} &  A_{\alpha'} \ar[d]^{f_{{\alpha},{\alpha'}}} \ar[r]^{h_{\alpha'}}&  X_{\alpha'}=X \ar@{=}[d] \ar[r]&0  \\
\bE_{\alpha}:\quad 0\ar[r] & B_{\alpha}\ar[r]^{k_{\alpha}}     &  A_{\alpha}  \ar[r]^{h_{\alpha}}
&  X_{\alpha}=X\ar[r]                              &0}$$
such that both $f_{{\alpha},{\alpha'}}$ and $g_{{\alpha},{\alpha'}}$ are epimorphisms.
\end{prt}

Set $\bE_0= 0\to 0\to 0\to 0\to 0$ and $\bE_1 = \bE$.
Suppose that $\alpha+1$ is a successor ordinal and we have $\bE_{\alpha}= \,0\to B_\alpha \to A_\alpha\to X\to0$.
Next we construct $\bE_{\alpha+1}$ in the following steps.

{\bf Step 1.} Construct $A_{\alpha+1}$ and $B_{\alpha+1}$.

Let $i\in V_{\alpha+1}\backslash V_\alpha$.
Then $\oplus_{a\in Q_{1}^{\ast\to i}}A_\alpha(s(a))\in\sf{X}$,
as each arrow $a\in Q_{1}^{\ast\to i}$ has the source $s(a)$ in $V_\alpha$;
see \corcite[2.8]{HJ19}.
Since the cotorsion pair $(\sf X, \sf Y)$ is complete,
there exists a short exact sequence
\begin{center}
$0\to \oplus_{a\in Q_{1}^{\ast\to i}}A_\alpha(s(a))\xra{(\epsilon_a)_{a\in Q_{1}^{\ast\to i}}} D_{\alpha+1}^i\to \overline{B}^i\to0$.
\end{center}
in $\A$ with $D_{\alpha+1}^i\in \sf X\cap\sf Y$ and $\overline{B}^i\in\sf{X}$.
For each arrow $a\in Q_{1}^{\ast\to i}$,
there exists the canonical inclusion
\begin{center}
$(\epsilon_a)_{a\in Q_{1}^{\ast\to i}}\circ \iota^{s(a),i}: \,\,
A_\alpha(s(a))\overset{\iota^{s(a),i}}{\hookrightarrow} \oplus_{a\in Q_{1}^{\ast\to i}}A_\alpha(s(a))
\xra{(\epsilon_a)_{a\in Q_{1}^{\ast\to i}}} D_{\alpha+1}^i$.
\end{center}
By assumption, one has $i\in Q_0\backslash V_\alpha$.
Hence $A_\alpha(a):A_\alpha(s(a))\to A_\alpha(i)=A(i)$.
So there is a canonical morphism
\begin{equation*}
\left(\begin{array}{cc}
A_\alpha(a) \\ (\epsilon_a)_{a\in Q_{1}^{\ast\to i}}\circ \iota^{s(a),i} \end{array}\right)
:A_\alpha(s(a))\to A(i)\oplus D_{\alpha+1}^i
\end{equation*}
which induces a monomorphism
\begin{equation}
\tag{$\dag$}
\left(\begin{array}{cc}(A_\alpha(a))_{a\in Q_{1}^{\ast\to i}} \\ (\epsilon_a)_{a\in Q_{1}^{\ast\to i}} \end{array}\right)
:\oplus_{a\in Q_{1}^{\ast\to i}}A_\alpha(s(a))\to A(i)\oplus D_{\alpha+1}^i,
\end{equation}
as $(\epsilon_a)_{a\in Q_{1}^{\ast\to i}}$ is a monomorphism.
Define $A_{\alpha+1}$ and $B_{\alpha+1}$ as follows:
\begin{equation*}
  A_{\alpha+1}(i)=
  \begin{cases}
    A_{\alpha}(i)                & \text{if $i\notin V_{\alpha+1}\backslash V_\alpha$}\,,\\
    \\
    A(i)\oplus D_{\alpha+1}^i    & \text{otherwise}.
  \end{cases}
\end{equation*}
and
\begin{equation*}
  B_{\alpha+1}(i)=
  \begin{cases}
    B_{\alpha}(i)                & \text{if $i\notin V_{\alpha+1}\backslash V_\alpha$}\,,\\
    \\
    B(i)\oplus D_{\alpha+1}^i    & \text{otherwise}.
  \end{cases}
\end{equation*}
For an arrow $a: j\to k$ in $Q_1$, the morphisms $A_{\alpha+1}(a)$ and $B_{\alpha+1}(a)$ are defined as:
\begin{rqm}
\item[--] If $k\in V_{\alpha+1}\backslash V_\alpha$, then $j\in V_\alpha$.
          Define
          $$A_{\alpha+1}(a)= \left(\begin{array}{cc}A_\alpha(a) \\ (\epsilon_a)_{a\in Q_{1}^{\ast\to k}}\circ \iota^{j,k}
          \end{array}\right): A_{\alpha+1}(j)=A_{\alpha}(j)\to A_{\alpha+1}(k)= A(k)\oplus D_{\alpha+1}^k,$$
          and define
          $$B_{\alpha+1}(a)= \left(\begin{array}{cc}B_\alpha(a) \\ (\epsilon_a)_{a\in Q_{1}^{\ast\to k}}\circ \iota^{j,k}\circ k_\alpha(j)
          \end{array}\right):B_{\alpha+1}(j)=B_{\alpha}(j)\to B_{\alpha+1}(k)= B(k)\oplus D_{\alpha+1}^k.$$
\item[--] If $j\in V_{\alpha+1}\backslash V_\alpha$, then $k\notin V_{\alpha+1}$; see \corcite[2.8]{HJ19}. Hence $A_{\alpha+1}(j)= A(j)\oplus D_{\alpha+1}^j$,
          $A_{\alpha+1}(k)= A_{\alpha}(k)=A(k)$,
          $B_{\alpha+1}(j)= B(j)\oplus D_{\alpha+1}^j$ and
          $B_{\alpha+1}(k)= B_{\alpha}(k)=B(k)$.
          Define $A_{\alpha+1}(a)$ as the composition of the projection $A(j)\oplus D_{\alpha+1}^j\twoheadrightarrow A(j)$
          followed by the morphism $A(j)\xra{A(a)}A(k)$, and define
          $B_{\alpha+1}(a)$ as the composition of the projection $B(j)\oplus D_{\alpha+1}^j\twoheadrightarrow B(j)$
          followed by the morphism $B(j)\xra{B(a)}B(k)$.

\item[--] For the other cases, define $A_{\alpha+1}(a)=A_{\alpha}(a)$ and $B_{\alpha+1}(a)=B_{\alpha}(a)$.
\end{rqm}
Note that $D_{\alpha+1}^i\in \sf X\cap\sf Y\subseteq\sf Y$
for each $i\in V_{\alpha+1}\backslash V_\alpha$.
It follows that $B(i)\oplus D_{\alpha+1}^i\in\sf{Y}$.
Thus $B_{\alpha+1}(i)\in \sf{Y}$ for all $i\in Q_0$.
On the other hand, for each $i\in V_{\alpha+1}\backslash V_\alpha$,
$A(i)\oplus D_{\alpha+1}^i\in \sf{X}$.
Hence $A_{\alpha+1}(i)\in \sf{X}$ for all $i\in Q_0$.
Thus, to show $A_{\alpha+1}$ and $B_{\alpha+1}$ constructed above
satisfy the desired conditions (b) and (c),
it remains to prove that $\C_i(A_{\alpha+1})\in \sf{X}$
for all $i\in V_{\alpha+1}\backslash V_\alpha$,
that is, the cokernel of the monomorphism
$$\varphi^{A_{\alpha+1}}_i=
\left(\begin{array}{cc}(A_\alpha(a))_{a\in Q_{1}^{\ast\to i}}\\ (\epsilon_a)_{a\in Q_{1}^{\ast\to i}}\end{array}\right)$$
given in $(\dag)$ is in $\sf{X}$.
Consider the following commutative diagram with exact rows and columns
$$\xymatrix@C=18pt@R=15pt{
&& A(i)\ar@{>->}[d]^{\iota}  \\
0 \ar[r] & \oplus_{a\in Q_{1}^{\ast\to i}}A_\alpha(s(a)) \ar[r]^{\ \ \ \ \varphi^{A_{\alpha+1}}_i}\ar@{=}[d] & A(i)\oplus D_{\alpha+1}^i \ar[r]\ar@{>>}[d]^{\pi} & \C_i(A_{\alpha+1})     \ar[r]\ar@.[d] & 0 \\
0 \ar[r] & \oplus_{a\in Q_{1}^{\ast\to i}}A_\alpha(s(a)) \ar[r]^{} &  D_{\alpha+1}^i  \ar[r] & \overline{B}^i     \ar[r]& 0. }
$$
Then there exists a morphism $\C_i(A_{\alpha+1})\to \overline{B}^i$
such that the right square is commutative.
By the Snake Lemma, one gets a short exact sequence
$0\to A(i)\to \C_i(A_{\alpha+1})\to \overline{B}^i\to 0$ in $\sf A$.
Since both $A(i)$ and $\overline{B}^i$ are in $\sf X$,
one concludes that $\C_i(A_{\alpha+1})\in \sf X$ as well.

{\bf Step 2.} Construct $\bE_{\alpha+1}$.

Set $X_{\alpha+1}=X$. Define $h_{\alpha+1}: A_{\alpha+1}\to X$ and $k_{\alpha+1}: B_{\alpha+1}\to A_{\alpha+1}$ as follows:
\begin{equation*}
  h_{\alpha+1}(i)=
  \begin{cases}
    h_{\alpha}(i)               & \text{if $i\notin V_{\alpha+1}\backslash V_\alpha$}\\
   \\
    (h_{\alpha}(i),0)    & \text{otherwise}
  \end{cases}
k_{\alpha+1}(i)=
  \begin{cases}
    k_{\alpha}(i)               & \text{if $i\notin V_{\alpha+1}\backslash V_\alpha$}\\

   \left( \begin{array}{cc}
                                     k_{\alpha}(i) & 0 \\
                                                      0 & 1 \\
                                                    \end{array}
                                                  \right) & \text{otherwise}.
  \end{cases}
\end{equation*}
According to the above constructions, one gets the desired exact sequence
$$\bE_{\alpha+1}:\,\,\,
0\to B_{\alpha+1}\xra{k_{\alpha+1}} A_{\alpha+1}\xra{h_{\alpha+1}} X_{\alpha+1}=X\to0.$$

{\bf Step 3.} Construct epimorphisms $f_{\alpha,{\alpha+1}}: A_{\alpha+1} \to A_{\alpha}$ and
$g_{\alpha,{\alpha+1}}: B_{\alpha+1} \to B_{\alpha}$.

Define $f_{\alpha,{\alpha+1}}$ and $g_{\alpha,{\alpha+1}}$ as follows:
\begin{equation*}
  f_{\alpha,{\alpha+1}}(i)=
  \begin{cases}
    \textrm{id}_{A_\alpha(i)}              & \text{if $i\notin V_{\alpha+1}\backslash V_\alpha$}\\
\\
    \left(1,  0 \right)    & \text{otherwise}
  \end{cases}
g_{\alpha,{\alpha+1}}(i)=
  \begin{cases}
    \textrm{id}_{B_\alpha(i)}              & \text{if $i\notin V_{\alpha+1}\backslash V_\alpha$}\\
\\
    \left( 1, 0\right)    & \text{otherwise}
  \end{cases}
\end{equation*}
It is clear that $f_{\alpha,{\alpha+1}}$ and $g_{\alpha,{\alpha+1}}$ are epimorphisms,
and it is easy to check that the diagram
$$\xymatrix{
\bE_{\alpha+1}:\quad 0\ar[r] &B_{\alpha+1}\ar[r]^{\quad k_{\alpha+1}\ \ \ }\ar[d]_{{g_{\alpha,{\alpha+1}}}} & A_{\alpha+1} \ar[d]_{f_{\alpha,{\alpha+1}}} \ar[r]^{h_{\alpha+1}}
& X \ar@{=}[d] \ar[r]&0  \\
\bE_{\alpha}:\quad \quad0\ar[r] &B_{\alpha}\ar[r]^{\quad k_{\alpha}\ \ \ \ \ }&  A_{\alpha}  \ar[r]^{h_{\alpha}} & X\ar[r]&0}$$
is commutative.

Suppose now that $\beta\leqslant\lambda$ is a limit ordinal and $\bE_\alpha$ is constructed for all $\alpha<\beta$.
Next we construct $\bE_\beta$.
In this case one has $V_\beta=\cup_{\alpha<\beta}V_{\alpha}$.
\begin{rqm}
\item[--] If $i\in V_\beta$, then $i\in V_\alpha$ for some ordinal $\alpha<\beta$,
          and so for each ordinal $\alpha<\alpha'\leqslant\beta$
          one has $\bE_{\alpha'}(i)=\bE_{\alpha}(i)$ as $i\notin V_{\alpha'}\backslash V_{\alpha}$.
\item[--] If $i\notin V_\beta$, then by hypothesis, one has $\bE_{\alpha}(i)=\bE(i)$ for all $\alpha<\beta$.
\end{rqm}
We let $\bE_\beta=\lim_{\alpha<\beta}\bE_\alpha$. Then one has
\begin{equation*}
  \bE_{\beta}(i)= \lim_{\alpha<\beta}\bE_\alpha(i)=
  \begin{cases}
    \bE_{\alpha}(i)\ \text{for some~} \alpha<\beta  & \text{if $i\in V_{\beta}$}\,,\\
    \bE(i) & \text{if $i\notin V_{\beta}$}\,.
  \end{cases}
\end{equation*}
Thus $\bE_{\beta}$ is an exact sequence in $\Rep{Q}{\A}$ satisfying the desired conditions.

Finally, let $A'=A_\lambda$ and $B'=B_\lambda$.
Then one gets the desired exact sequence $\bE'=\bE_\lambda:\,\,\,0\to B'\to A'\to X\to 0$.
\qed
\end{bfhpg}

The next result is immediate by Theorem \ref{AAA} and Remark \ref{equality}.

\begin{cor}\label{corA}
Let $\sf A$ be an abelian category, and let $(\sf X,\sf Y)$ be a complete cotorsion pair in $\sf A$.
\begin{prt}
\item Suppose that $Q$ is a left rooted quiver and $\sf A$ satisfies the axiom AB4. Then $(\Phi(\sf X),\Rep{Q}{\sf Y})$ is a complete cotorsion pair in the category $\Rep{Q}{\A}$.
\item Suppose that $Q$ is a right rooted quiver and $\sf A$ satisfies the axiom AB4*. Then $(\Rep{Q}{\sf X},\Psi(\sf Y))$ is a complete cotorsion pair in the category $\Rep{Q}{\A}$.
\end{prt}
\end{cor}

The following result shows that under some mild conditions the converses of the statements in Corollary \ref{corA} hold true.

\begin{prp}\label{If-part complete}
Let $\sf X$ be a subcategory of $\sf A$, which is closed under extensions.
\begin{prt}
\item Suppose that $Q$ is a left rooted quiver, $\sf X$ is closed under coproducts and $\sf A$ satisfies the axiom $AB4$.
      If $(\Phi({\sf X}),\Rep{Q}{\sf{X}^\perp})$ is a complete cotorsion pair in $\Rep{Q}{\A}$,
      then $(\sf X,\sf{X}^\perp)$ is a complete cotorsion pair in $\sf A$.
\item Suppose that $Q$ is a right rooted quiver, $\sf X$ is closed under products and $\sf A$ satisfies the axiom $AB4^*$.
      If $(\Rep{Q}{{^\perp\sf{X}}},\Psi({\sf X}))$ is a complete cotorsion pair in $\Rep{Q}{\A}$,
      then $({^\perp\sf{X}},\sf X)$ is a complete cotorsion pair in $\sf A$.
\end{prt}
\end{prp}

\begin{prf*}
We only prove (a); the statement (b) is proved dually.

By Proposition \ref{another side of cp}(a), the pair $(\sf X,\sf{X}^\perp)$ forms a cotorsion pair in $\sf A$. We next show that this cotorsion pair is complete.
To this end, let $M$ be an object in $\sf A$.
Consider the representation ${\sf s}_i(M)$ in $\Rep{Q}{\A}$
where $i$ is an fixed vertex in $Q_0$.
Since the cotorsion pair $(\Phi({\sf X}),\Rep{Q}{\sf{X}^\perp})$ is complete by assumption,
there exist two short exact sequences
\begin{center}
$0\to Y\to X\to {\sf s}_i(M)\to0$
\quad and \quad
$0\to{\sf s}_i(M)\to Y'\to X'\to0$
\end{center}
in $\Rep{Q}{\A}$ such that both $X$ and $X'$ are in $\Phi({\sf X})$ and
both $Y$ and $Y'$ are in $\Rep{Q}{\sf{X}^\perp}$.
Note that ${\sf s}_i(M)(i)=M$.
We have two short exact sequences
\begin{center}
$0\to Y(i)\to X(i)\to M\to0$
\quad and \quad
$0\to M\to Y'(i)\to X'(i)\to0$
\end{center}
in $\A$ such that both $Y(i)$ and $Y'(i)$ are in $\sf{X}^\perp$.
Moreover, since $Q$ is left rooted and $\sf X$ is closed under extensions and coproducts,
we conclude by Lemma \ref{X(i) in X}(a) that both $X(i)$ and $X'(i)$ are in $\sf{X}$.
Thus, the cotorsion pair $(\sf X,\sf{X}^\perp)$ is complete.
\end{prf*}

As an immediately consequence of Corollary \ref{corA} and Proposition \ref{If-part complete}, we have the next result.

\begin{cor}\label{if and only if C}
Let $\sf X$ and $\sf Y$ be two subcategories in $\sf A$.
\begin{prt}
\item Suppose that $\sf A$ satisfies the axiom $AB4$.
      If $Q$ is a left rooted quiver,
      then $(\sf X,\sf{X}^\perp)$ is a complete cotorsion pair in $\sf A$ if and only if
      $(\Phi({\sf X}),\Rep{Q}{\sf{X}^\perp})$ is a complete cotorsion pair in $\Rep{Q}{\A}$ and
      $\sf X$ is closed under extensions and coproducts.
\item Suppose that $\sf A$ satisfies the axiom $AB4^*$.
      If $Q$ is a right rooted quiver,
      then $({^\perp\sf{X}},\sf X)$ is a complete cotorsion pair in $\sf A$ if and only if
      $(\Rep{Q}{{^\perp\sf{X}}},\Psi({\sf X}))$ is a complete cotorsion pair in $\Rep{Q}{\A}$ and
      $\sf X$ is closed under extensions and products.
\end{prt}
\end{cor}

\section{Applications to special precovering/preenveloping subcategories}\label{model}
\noindent
In this section we aim at proving Theorem \ref{BBBB} in the introduction.
The next result is a key to do that, which is from Hu and Zhu \lemcite[3.9]{Hu-Zhu}.

For a subcategory $\sf L$ of $\A$,
we denote by $\Smd(\sf L)$ the subcategory of $\A$
consisting of all direct summands of objects in $\sf L$.

\begin{lem}\label{Key bridge}
Let $\sf L$ be a subcategory of $\A$.
\begin{prt}
\item Suppose that $\A$ has enough injectives.
      If $\sf L$ is special precovering in $\A$,
      then $(\Smd(\sf L),\sf L^\perp)$ forms a complete cotorsion pair.
\item Suppose that $\A$ has enough projectives.
      If $\sf L$ is special preenveloping in $\sf A$,
      then $({^\perp\sf L},\Smd(\sf L))$ forms a complete cotorsion pair.
\end{prt}
\end{lem}

\begin{rmk}
An interesting observation from the above result is that a special precovering (resp., preenveloping) subcategory,
which is closed under direct summands,
in an abelian category with enough injectives (resp., projectives)
must be closed under extensions.
\end{rmk}

With the help of Lemma \ref{Key bridge},
we give the proof of Theorem \ref{BBBB} in the introduction.

\begin{bfhpg}[\bf Proof of Theorem \ref{BBBB}]\label{PROOF OF B1}
We only prove (a); one can prove (b) dually.

Since $\sf L$ is a special precovering subcategory in $\sf A$ by assumption,
one concludes by Lemma \ref{Key bridge}(a) that
$(\Smd(\sf L),\sf L^\perp)$ is a complete cotorsion pair in $\A$.
Hence, $(\Phi(\Smd(\sf L)), \Rep{Q}{\sf L^\perp})$ is a complete cotorsion pair in $\Rep{Q}{\A}$
by Corollary \ref{corA}(a).
This implies that for any representation $M\in\Rep{Q}{\A}$,
there exists a short exact sequence
$$\bE:\,\,\,0\to K \overset{}\rightarrow N \overset{}\rightarrow M \to0$$
in $\Rep{Q}{\A}$ with $N\in \Phi(\Smd(\sf L))$ and $K\in \Rep{Q}{\sf L^\perp}$.
Note that $(\Smd(\sf L),\sf L^\perp)$ is a cotorsion pair. Then $\Smd(\sf L)$ is closed under extensions and coproducts, and so $N(i)\in \Smd(\sf L)$ for any vertex $i\in Q_0$ by Lemma \ref{X(i) in X}(a).
But $N$ may fail in $\Phi(\sf L)$.
We will use transfinite induction to construct a short exact sequence
$$\bE':\,\,\,0\to K'\to N'\to M\to 0$$
in $\Rep{Q}{\A}$ such that
$N'\in \Phi(\sf L)$ and $K'\in \Rep{Q}{\sf L^\perp}$.
If we have done, then $\Phi(\sf L)$ is special precovering in $\Rep{Q}{\A}$, as $\Rep{Q}{\sf L^\perp}= \Phi(\Smd(\sf L))^\perp\subseteq\Phi(\sf L)^\perp$.

Let $\{V_{\alpha}\}$ be the transfinite sequence of subsets of $Q_0$.
Since $Q$ is left rooted, one has $Q_0=V_{\lambda}$ for some ordinal $\lambda$.
Next we construct a continuous inverse $\lambda$-sequence
$$\{\,\bE_\alpha:\,\,\, 0\to K_\alpha \overset{g_\alpha}\longrightarrow N_\alpha\overset{f_\alpha}\longrightarrow M_\alpha\to0 \,\}_{\alpha\leqslant\lambda}$$
of short exact sequences in $\Rep{Q}{\A}$ satisfying the following conditions:

\begin{prt}
\item For every ordinal $0<\alpha\leqslant\lambda$, $M_\alpha=M$.
\item For every ordinal $0<\alpha\leqslant\lambda$ and for each $i\in Q_0\backslash V_\alpha$,
 $N_\alpha(i) = N(i)\in \Smd(\sf L)$ and $K_\alpha(i) = K(i)\in \sf L^\perp$.
\item For every ordinal $\alpha\leqslant\lambda$ and each $i\in V_\alpha$,
$\varphi_i^{N_\alpha}$ is a monomorphism, $\C_i(N_\alpha)\in \sf{L}$ and $K_\alpha(i)\in \sf L^\perp$.
\item For every $0<\alpha<\alpha'\leqslant\lambda$,
$N_\alpha(i) = N_{\alpha'}(i)$ and $K_\alpha(i) = K_{\alpha'}(i)$
for all $i\notin V_{\alpha'}\backslash V_\alpha$,
and there exists the following commutative diagram in $\Rep{Q}{\A}$
$$\xymatrix{
\bE_{\alpha'}:\quad    0\ar[r] & K_{\alpha'}\ar[r]^{g_{\alpha'}}\ar[d]^{{h_{\alpha,{\alpha'}}}} &  N_{\alpha'} \ar[d]^{k_{\alpha,{\alpha'}}} \ar[r]^{f_{\alpha'}\quad}&  M_{\alpha'}=M \ar@{=}[d] \ar[r]&0  \\
\bE_{\alpha}:\quad 0\ar[r] & K_{\alpha}\ar[r]^{g_{\alpha}}     &  N_{\alpha}  \ar[r]^{f_{\alpha}\quad}
&  M_{\alpha}=M\ar[r]                              &0}$$
such that both $h_{\alpha,{\alpha'}}$ and $k_{\alpha,{\alpha'}}$ are epimorphisms.
\end{prt}

Set $\bE_0= 0\to 0\to 0\to 0\to 0$ and $\bE_1 = \bE$.
Suppose that $\alpha+1$ is a successor ordinal and
we have $\bE_{\alpha}= \,0\to K_\alpha \to N_\alpha\to M\to0$.
Next we construct $\bE_{\alpha+1}$ in the following steps.

{\bf Step 1.} Construct $N_{\alpha+1}$ and $K_{\alpha+1}$.

Let $i\in V_{\alpha+1}\backslash V_\alpha$.
Since the cotorsion pair $(\Smd(\sf L),\sf L^\perp)$ is complete,
there exists a short exact sequence
\begin{center}
$0\to \oplus_{a\in Q_{1}^{\ast\to i}}N_\alpha(s(a))\xra{(\epsilon_a)_{a\in Q_{1}^{\ast\to i}}} D_{\alpha+1}^i\to \overline{B}^i\to0$.
\end{center}
in $\A$ with $D_{\alpha+1}^i\in\sf L^\perp$ and $\overline{B}^i\in\Smd(\sf L)$.
For each $a\in Q_{1}^{\ast\to i}$,
there exists the canonical inclusion
\begin{center}
$(\epsilon_a)_{a\in Q_{1}^{\ast\to i}}\circ \iota^{s(a),i}: \,\,
N_\alpha(s(a))\overset{\iota^{s(a),i}}{\hookrightarrow} \oplus_{a\in Q_{1}^{\ast\to i}}N_\alpha(s(a))
\xra{(\epsilon_a)_{a\in Q_{1}^{\ast\to i}}} D_{\alpha+1}^i$.
\end{center}
By assumption, one has $i\in Q_0\backslash V_\alpha$.
Hence, $N_\alpha(a):N_\alpha(s(a))\to N_\alpha(i)=N(i)$.
So there is a canonical morphism
\begin{equation*}
\left(\begin{array}{cc}
N_\alpha(a) \\ (\epsilon_a)_{a\in Q_{1}^{\ast\to i}}\circ \iota^{s(a),i} \end{array}\right)
:N_\alpha(s(a))\to N(i)\oplus D_{\alpha+1}^i
\end{equation*}
which induces a monomorphism
\begin{equation}
\tag{$\dag$}
\varphi=\left(\begin{array}{cc}(N_\alpha(a))_{a\in Q_{1}^{\ast\to i}} \\ (\epsilon_a)_{a\in Q_{1}^{\ast\to i}} \end{array}\right)
:\oplus_{a\in Q_{1}^{\ast\to i}}N_\alpha(s(a))\to N(i)\oplus D_{\alpha+1}^i.
\end{equation}

Consider the following commutative diagram with exact rows and columns
$$\xymatrix@C=18pt@R=15pt{
&& N(i)\ar@{>->}[d]^{\iota}  \\
0 \ar[r] & \oplus_{a\in Q_{1}^{\ast\to i}}N_\alpha(s(a)) \ar[r]^{\ \ \ \ \varphi}\ar@{=}[d] & N(i)\oplus D_{\alpha+1}^i \ar[r]\ar@{>>}[d]^{\pi} & \coker(\varphi)     \ar[r]\ar@.[d] & 0 \\
0 \ar[r] & \oplus_{a\in Q_{1}^{\ast\to i}}N_\alpha(s(a)) \ar[r]^{} &  D_{\alpha+1}^i  \ar[r] & \overline{B}^i     \ar[r]& 0. }
$$
Then there exists a morphism $\coker(\varphi)\to \overline{B}^i$
such that the right square is commutative.
By the Snake Lemma, one gets a short exact sequence
$0\to N(i)\to \coker(\varphi)\to \overline{B}^i\to 0$ in $\sf A$.
Since both $N(i)$ and $\overline{B}^i$ are in $\Smd(\sf L)$,
we conclude that $\coker(\varphi)\in \Smd(\sf L)$ as well.
Moreover, since $\sf L$ is special precovering in $\sf A$ by assumption,
there exists a short exact sequence
\begin{center}
$0\to H_{\alpha+1}^i \to L_{\alpha+1}^i\to \coker(\varphi)\to0$.
\end{center}
in $\A$ with $L_{\alpha+1}^i\in\sf L$ and $H_{\alpha+1}^i\in\sf L^\perp = \Smd(\sf L)^\perp$.
Clearly, this short exact sequence is split.
Hence, one has
\begin{equation}
\tag{$\ddag$}
\coker(\varphi)\oplus H_{\alpha+1}^i \cong L_{\alpha+1}^i \in \sf L.
\end{equation}

Define $N_{\alpha+1}$ and $K_{\alpha+1}$ as follows:
\begin{equation*}
  N_{\alpha+1}(i)=
  \begin{cases}
    N_{\alpha}(i)                                     & \text{if $i\notin V_{\alpha+1}\backslash V_\alpha$}\,,\smallskip\\
    N(i)\oplus D_{\alpha+1}^i \oplus H_{\alpha+1}^i   & \text{otherwise}.
  \end{cases}
\end{equation*}
and
\begin{equation*}
  K_{\alpha+1}(i)=
  \begin{cases}
    K_{\alpha}(i)                                     & \text{if $i\notin V_{\alpha+1}\backslash V_\alpha$}\,,\smallskip\\
    K(i)\oplus D_{\alpha+1}^i \oplus H_{\alpha+1}^i   & \text{otherwise}.
  \end{cases}
\end{equation*}
For an arrow $a: j\to k$ in $Q_1$, the morphisms $N_{\alpha+1}(a)$ and $K_{\alpha+1}(a)$ are defined as:
\begin{rqm}
\item[--] If $k\in V_{\alpha+1}\backslash V_\alpha$, then $j\in V_\alpha$; see \corcite[2.8]{HJ19}.
          Define
          $$N_{\alpha+1}(a)= \left(\begin{array}{cc}N_\alpha(a) \smallskip\\  (\epsilon_a)_{a\in Q_{1}^{\ast\to k}}\circ \iota^{j,k}\smallskip \\ 0
          \end{array}\right) :$$
          $$N_{\alpha+1}(j)=N_{\alpha}(j)\longrightarrow N_{\alpha+1}(k)= N(k)\oplus D_{\alpha+1}^k\oplus H_{\alpha+1}^i,$$
          and define
          $$K_{\alpha+1}(a)= \left(\begin{array}{cc}K_\alpha(a) \\ (\epsilon_a)_{a\in Q_{1}^{\ast\to k}}\circ \iota^{j,k}\circ g_\alpha(j) \smallskip \\ 0
          \end{array}\right) :$$
          $$K_{\alpha+1}(j)=K_{\alpha}(j)\longrightarrow K_{\alpha+1}(k)= K(k)\oplus D_{\alpha+1}^k\oplus H_{\alpha+1}^i.$$
\item[--] If $j\in V_{\alpha+1}\backslash V_\alpha$, then $k\notin V_{\alpha+1}$.
          Hence $N_{\alpha+1}(j)= N(j)\oplus D_{\alpha+1}^j\oplus H_{\alpha+1}^j$,
          $N_{\alpha+1}(k)= N_{\alpha}(k)=N(k)$,
          $K_{\alpha+1}(j)= K(j)\oplus D_{\alpha+1}^j\oplus H_{\alpha+1}^j$ and
          $K_{\alpha+1}(k)= K_{\alpha}(k)=K(k)$.
          Define $N_{\alpha+1}(a)$ as the composition of the projection
          $N(j)\oplus D_{\alpha+1}^j\oplus H_{\alpha+1}^j\twoheadrightarrow N(j)$ followed by $N(j)\xra{N(a)}N(k)$, and define
          $K_{\alpha+1}(a)$ as the composition of the projection
          $K(j)\oplus D_{\alpha+1}^j\oplus H_{\alpha+1}^j\twoheadrightarrow K(j)$ followed by $K(j)\xra{K(a)}K(k)$.

\item[--] For the other cases, define $N_{\alpha+1}(a)=N_{\alpha}(a)$ and $K_{\alpha+1}(a)=K_{\alpha}(a)$.
\end{rqm}

Note that for each $i\in V_{\alpha+1}\backslash V_\alpha$, both $D_{\alpha+1}^i$ and $H_{\alpha+1}^i$ are in $\sf L^\perp$.
It follows that $K(i)\oplus D_{\alpha+1}^i \oplus H_{\alpha+1}^i\in\sf L^\perp$.
Thus $K_{\alpha+1}(i)\in \sf L^\perp$ for all $i\in Q_0$.
On the other hand, for $i\in V_{\alpha+1}\backslash V_\alpha$, since $(\dag)$ is a monomorphism,
we see that
$$\varphi^{N_{\alpha+1}}_i=
\left(\begin{array}{cc}(N_\alpha(a))_{a\in Q_{1}^{\ast\to i}}\\ (\epsilon_a)_{a\in Q_{1}^{\ast\to i}}\smallskip\\ 0\end{array}\right):
\oplus_{a\in Q_{1}^{\ast\to i}}N_\alpha(s(a))\to N(i)\oplus D_{\alpha+1}^i\oplus H_{\alpha+1}^i
$$
is a monomorphism as well.
Thus, to show $N_{\alpha+1}$ and $K_{\alpha+1}$ constructed above
satisfy the desired conditions (b) and (c),
it remains to prove that $\C_i(N_{\alpha+1})\in \sf{L}$
for all $i\in V_{\alpha+1}\backslash V_\alpha$.
Indeed, according to what we constructed above,
it is clear that $\C_i(N_{\alpha+1})\cong \coker(\varphi)\oplus H_{\alpha+1}^i \cong L_{\alpha+1}^i\in \sf L$ by $(\ddagger)$,
as desired.

{\bf Step 2.} Construct $\bE_{\alpha+1}$.

Set $M_{\alpha+1}=M$. Define $f_{\alpha+1}: N_{\alpha+1}\to M$ and $g_{\alpha+1}: K_{\alpha+1}\to N_{\alpha+1}$ as follows:
\begin{equation*}
  f_{\alpha+1}(i)=
  \begin{cases}
    f_{\alpha}(i)               & \text{if $i\notin V_{\alpha+1}\backslash V_\alpha$}\,,\\
   \\
    \left(f_{\alpha}(i), 0, 0\right)    & \text{otherwise}.
  \end{cases}
\end{equation*}
and
\begin{equation*}
 \quad \, g_{\alpha+1}(i)=
  \begin{cases}
    g_{\alpha}(i)               & \text{if $i\notin V_{\alpha+1}\backslash V_\alpha$}\,,\\
   \\
   \left( \begin{array}{cc}
                                              g_{\alpha}(i)  & 0 \quad 0 \\
                                              0              & 1 \quad 0 \\
                                              0              & 0 \quad 1
                                                    \end{array}
                                                  \right) & \text{otherwise}.
  \end{cases}
\end{equation*}
According to the above constructions, one gets the desired exact sequence
$$\bE_{\alpha+1}:\,\,\,
0\to K_{\alpha+1}\xra{g_{\alpha+1}} N_{\alpha+1}\xra{f_{\alpha+1}} M_{\alpha+1}=M\to0.$$

{\bf Step 3.}
Construct epimorphisms $k_{\alpha,{\alpha+1}}: N_{\alpha+1}\to N_\alpha $ and $h_{\alpha,{\alpha+1}}: K_{\alpha+1} \to K_\alpha$.

Define $k_{\alpha,{\alpha+1}}$ and $h_{\alpha,{\alpha+1}}$ as follows:
\begin{equation*}
  k_{\alpha,{\alpha+1}}(i)=
  \begin{cases}
    \textrm{id}_{N_\alpha(i)}              & \text{if $i\notin V_{\alpha+1}\backslash V_\alpha$}\,\\
\\
    \left(1, 0, 0 \right)    & \text{otherwise}
  \end{cases}
h_{\alpha,{\alpha+1}}(i)=
  \begin{cases}
    \textrm{id}_{K_\alpha(i)}              & \text{if $i\notin V_{\alpha+1}\backslash V_\alpha$}\,\\
\\
    \left(1, 0, 0 \right)    & \text{otherwise}
  \end{cases}
\end{equation*}
It is clear that $k_{\alpha,{\alpha+1}}$ and $h_{\alpha,{\alpha+1}}$ are epimorphisms,
and it is easy to check that the diagram
$$\xymatrix{
\bE_{\alpha+1}:\quad0\ar[r] &K_{\alpha+1}\ar[r]^{\quad g_{\alpha+1}\ \ \ }\ar[d]_{{h_{{\alpha},{\alpha+1}}}} & N_{\alpha+1} \ar[d]_{k_{{\alpha},{\alpha+1}}} \ar[r]^{f_{\alpha+1}}
& M \ar@{=}[d] \ar[r]&0  \\
\bE_{\alpha}:\quad\quad 0\ar[r] &K_{\alpha}\ar[r]^{\quad g_{\alpha}\ \ \  }&  N_{\alpha}  \ar[r]^{f_{\alpha}} & M\ar[r]&0}$$
is commutative.

Suppose now that $\beta\leqslant\lambda$ is a limit ordinal and $\bE_\alpha$ is constructed for all $\alpha<\beta$.
Next we construct $\bE_\beta$.
In this case one has $V_\beta=\cup_{\alpha<\beta}V_{\alpha}$.
\begin{rqm}
\item[--] If $i\in V_\beta$, then $i\in V_\alpha$ for some ordinal $\alpha<\beta$,
          and so for each ordinal $\alpha<\alpha'\leqslant\beta$
          one has $\bE_{\alpha'}(i)=\bE_{\alpha}(i)$ as $i\notin V_{\alpha'}\backslash V_{\alpha}$.
\item[--] If $i\notin V_\beta$, then by hypothesis, one has $\bE_{\alpha}(i)=\bE(i)$ for all $\alpha<\beta$.
\end{rqm}
We let $\bE_\beta=\lim_{\alpha<\beta}\bE_\alpha$. Then one has
\begin{equation*}
  \bE_{\beta}(i)= \lim_{\alpha<\beta}\bE_\alpha(i)=
  \begin{cases}
    \bE_{\alpha}(i)\ \text{for some~} \alpha<\beta  & \text{if $i\in V_{\beta}$}\,,\\
    \bE(i) & \text{if $i\notin V_{\beta}$}\,.
  \end{cases}
\end{equation*}
Thus $\bE_{\beta}$ is a short exact sequence in $\Rep{Q}{\A}$ satisfying the desired conditions.

Finally, let $A'=A_\lambda$ and $B'=B_\lambda$.
Then one gets the desired short exact sequence $\bE'=\bE_\lambda:\,\,\,0\to B'\to A'\to X\to 0$.
\qed
\end{bfhpg}

\begin{rmk}
We don't assume that the subcategory $\sf L$ is closed under direct summands in Theorem \ref{BBBB}. Actually, if $\sf L$ is closed under direct summands, then one has $\Smd(\sf L)=\sf L$, and so under the assumptions of Theorem \ref{BBBB} the pair $(\sf L,\sf L^\perp)$ forms a complete cotorsion pair; see Lemma \ref{Key bridge}. Now Theorem \ref{AAA} yields that $(\Phi(\sf L),\Phi(\sf L)^\perp)$ is a complete cotorsion pair in $\Rep{Q}{\A}$. Hence $\Phi(\sf L)$ is special precovering.

We notice that there do exist a large number of special precovering subcategories
that may not be closed under direct summands. For example, the subcategory of $\Mod{R}$ consisting of all free $R$-modules is so, where $R$ is an associated ring and $\Mod{R}$ denote the category of all $R$-modules.
\end{rmk}

We end this section with the next result, which asserts that under some mild conditions the converses of the statements in Theorem \ref{BBBB} hold true.

\begin{prp}\label{an-side of sp}
Let $\sf L$ be a subcategory of $\sf A$, which is closed under extensions.
\begin{prt}
\item Suppose that $Q$ is a left rooted quiver and
      $\sf A$ satisfies the axiom $AB4$.
      If $\Phi({\sf L})$ is special precovering in $\Rep{Q}{\A}$
      and $\sf L$ is closed under coproducts,
      then $\sf L$ is special precovering in $\sf A$.
\item Suppose that $Q$ is a right rooted quiver and
      $\sf A$ satisfies the axiom $AB4^*$.
      If $\Psi({\sf L})$ is special preenveloping in $\Rep{Q}{\A}$
      and $\sf L$ is closed under products,
      then $\sf L$ is special preenveloping in $\sf A$.
\end{prt}
\end{prp}
\begin{prf*}
We only prove (a); one can prove (b) dually.

We claim that $\Phi({\sf L})^\perp\subseteq\Rep{Q}{{\sf L}^\perp}$.
To this end, let $X$ be a representation in $\Phi({\sf L})^\perp$ and $L$ an object in $\sf L$.
By \prpcite[5.2(a)]{HJ19},
we see that
$$\textrm{Ext}_{\A}^1(L,{\sf e}_i(X))\cong\textrm{Ext}_{\Rep{Q}{\A}}^1({\sf f}_i(L),X)$$
for any vertex $i\in Q_0$.
Since $\sf L$ is closed under coproducts by assumption,
it is easy to see that ${\sf f}_i(L)\in\Phi({\sf L})$.
Hence, $\textrm{Ext}_{\Rep{Q}{\A}}^1({\sf f}_i(L),X)=0$,
which implies that $\textrm{Ext}_{\A}^1(L,{\sf e}_i(X))=0$.
Note that ${\sf e}_i(X)=X(i)$.
It follows that $X(i)\in {\sf L}^\perp$.
Thus, $X\in \Rep{Q}{{\sf L}^\perp}$.
This yields that $\Phi({\sf L})^\perp\subseteq\Rep{Q}{{\sf L}^\perp}$.

Now let $M$ be an object in $\sf A$.
Since $\Phi({\sf L})$ is special precovering in $\Rep{Q}{\A}$,
there exists for an fixed vertex $i\in Q_0$ a short exact sequence
$$0\to K \to H \to {\sf s}_i(M)\to0$$
in $\Rep{Q}{\A}$ with $H\in \Phi({\sf L})$ and $K\in \Phi({\sf L})^\perp\subseteq\Rep{Q}{{\sf L}^\perp}$.
Note that $Q$ is left rooted and $\sf L$ is closed under extensions and coproducts by assumptions.
It follows from Lemma \ref{X(i) in X}(a) that $H(i)\in \sf L$.
Hence, we obtain a short exact sequence
$$0\to K(i) \to H(i) \to {\sf s}_i(M)(i)=M\to0$$
such that $H(i)\in \sf L$ and $K(i)\in \sf L^\perp$.
This yields that $\sf L$ is special precovering.
\end{prf*}

\appendix
\section*{Appendix}
\stepcounter{section}
\noindent
In this section we reprove Theorem \ref{AAA}(a) over a simple left rooted quiver
to comprehend the idea of the proof.

Let $Q$ be the quiver
$$
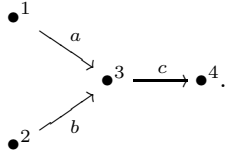
\xymatrix@C=20pt@R=10pt{
  \bullet^{1} \ar[dr]^{a} \\
  & \bullet^{3} \ar[r]^{c}  &  \bullet^{4}. \\
  \bullet^{2}  \ar[ur]_{b} }$$
Then $Q$ is a left rooted quiver with $V_0=\emptyset$, $V_1=\{1,2\}$, $V_2=\{1,2,3\}$, $V_3=Q_0$.

\begin{bfhpg*}[\bf Proof of Theorem A(a)]\label{pf1}
Let $X$ be a representation in $\Rep{Q}{\A}$.
Since $(\sf X,\sf Y)$ is a complete cotorsion pair in $\sf A$,
one concludes that for each vertex $i\in \{1,2,3,4\}$
there exists a short exact sequence
$$0\to B(i) \overset{k(i)}\longrightarrow A(i) \overset{h(i)}\longrightarrow X(i) \to 0$$
in $\A$ with $A(i)\in \sf X$ and $B(i)\in\sf Y$.
Moreover, for each arrow $e: i\to j \in \{a,b,c\}$ and for each morphism
$\xymatrix{A(i)\ar[rr]^{ X(e)\circ h(i)}& & X(j)}$
there exists a morphism $A(e):A(i)\to A(j)$ (and hence a morphism $B(e):B(i)\to B(j)$)
such that the diagram
$$\xymatrix{
0\ar[r] &B(i)\ar[r]^{k(i)}\ar@.[d]_{B(e)}& A(i) \ar@.[d]_{A(e)} \ar[r]^{h(i)}&  X(i) \ar[d]^{X(e)}\ar[r]&0  \\
0\ar[r] &B(j)\ar[r]^{k(j)}&  A(j)  \ar[r]^{h(j)} & X(j)\ar[r]&0}$$
is commutative.
Hence one gets an exact sequence
$$\bE:\,\,\,0\to B\overset{k}\longrightarrow A\overset{h}\longrightarrow X\to 0$$
in $\Rep{Q}{\A}$ such that $A\in \Rep{Q}{\sf X}$ and $B\in \Rep{Q}{\sf Y}$, which is in the following form:
$$\xymatrix@R=0.005cm{
  &&& 0 \ar[ddd] \\
  &0\ar[ddd]&&&& 0 \ar[ddd] &  0 \ar[ddd]\\
  &&&& 0 \ar[ddd]\\
  &&& B(1)\ar@.[drr]^{\quad \quad B(a)}\ar[ddd]_{k(1)} \\
  &~B:\ar[ddd]_{k}&&&& B(3) \ar@.[r]^{B(c)} \ar[ddd]_{k(3)} & B(4)\ar[ddd]_{k(4)} \\
  &&&& B(2)  \ar@.[ur]_{X(b)} \ar[ddd]_{k(2)}\\
  &&& A(1) \ar@.[drr]^{\quad \quad A(a)} \ar[ddd]_{h(1)}\\
  &~A:\ar[ddd]_{h}&&&& A(3) \ar@.[r]^{A(c)} \ar[ddd]_{h(3)} &  A(4)\ar[ddd]_{h(4)} \\
  &&&& A(2)  \ar@.[ur]_{A(b)}\ar[ddd]_{h(2)}\\
  &&& X(1) \ar[drr]^{\quad \quad X(a)} \ar[ddd]\\
  &~X:\ar[ddd]&&&& X(3) \ar[r]^{X(c)} \ar[ddd] &  X(4)\ar[ddd]\\
  &&&& X(2)  \ar[ur]_{X(b)}\ar[ddd]\\
  &&&0  \\
  &0&&&& 0  &  0 \\
  &&&&0
  }$$
We let $\bE_1=\bE$. It is clear that $\varphi^A_1$ and $\varphi^A_2$ are monomorphisms
with both $\C_1(A)=A(1)$ and $\C_2(A)=A(2)$ in $\sf X$,
as the vertices 1 and 2 are in $V_1$.
But $\varphi^A_3$ and $\varphi^A_4$ may fail to be monomorphisms,
and $\C_3(A)$ and $\C_4(A)$ may fail to be in $\sf X$.
Next we repair $\bE_1$ at the vertex 3.

Since $A(1)\oplus A(2)\in\sf X$ and the cotorsion pair $(\sf X, \sf{Y})$ is complete,
there exists a short exact sequence
\begin{center}
$0\to A(1)\oplus A(2)\xra{(\epsilon_2^1,\epsilon_2^2)} D_{2}^3\to {\overline{B}^3}\to0$.
\end{center}
in $\A$ with $D_{2}^3\in \sf{X}\cap\sf{Y}$ and ${\overline{B}}^3\in{\sf{X}}$.
Let $B_2$ and $A_2$ be representations as follows:
$$\xymatrix@R=0.05cm{
  &&B(1) \ar[dr]^{B_2(a)} \\
   &B_2:&& B(3)\oplus D_2^3 \ar[r]^{\quad B_2(c)}  &  B(4) \\
  &&B(2)  \ar[ur]_{B_2(b)} }$$
and
$$\xymatrix@R=0.05cm{
  &&A(1) \ar[dr]^{A_2(a)} \\
  &A_2:&& A(3)\oplus D_2^3 \ar[r]^{\quad A_2(c)}  &  A(4), \\
  &&A(2)  \ar[ur]_{A_2(b)} }$$
where $B_2(a)= \left(\begin{array}{cc}B(a) \\
                                    \epsilon_2^1\circ k(1)
                       \end{array}\right)$,
$B_2(b)= \left(\begin{array}{cc}B(b) \\
                                    \epsilon_2^2\circ k(2)
                       \end{array}\right)$,
$B_2(c)=\left(B(c),0\right)$,
$A_2(a)= \left(\begin{array}{cc}A(a) \\
                                    \epsilon_2^1
                       \end{array}\right)$,
$A_2(b)= \left(\begin{array}{cc}A(b) \\
                                    \epsilon_2^2
                       \end{array}\right)$,
and $A_2(c)=\left(A(c),0\right)$.

Note that $D_{2}^3\in{\sf{Y}}$.
It follows that $B(3)\oplus D_{2}^3\in\sf{Y}$, and so one has $B_2\in \Rep{Q}{\sf{Y}}$.
On the other hand, the morphism
$$\varphi^{A_2}_3=
\left(\begin{array}{cc}
A(a) \quad A(b)\\
\epsilon_2^1 \quad \epsilon_2^2
      \end{array}\right)$$
is now a monomorphism, as $(\epsilon_2^1,\epsilon_2^2)$ is so.
Next, we prove that $\C_3(A_{2})= \Coker\,\varphi^{A_2}_3$ is in ${\sf{X}}$.
To this end, consider the following commutative diagram with exact rows and columns:
$$\xymatrix@C=15pt@R=15pt{
&& A(3)\ar@{>->}[d]^{\iota}  \\
0 \ar[r] & A(1)\oplus A(2) \ar[r]^{\ \ \varphi^{A_2}_3}\ar@{=}[d] & A(3)\oplus D_{2}^3 \ar[r]\ar@{>>}[d]^{\pi} & \C_3(A_{2})     \ar[r]\ar@.[d] & 0 \\
0 \ar[r] & A(1)\oplus A(2) \ar[r]^{} &  D_{2}^3  \ar[r] & \overline{B}^3 \ar[r]& 0. }
$$
By the Snake Lemma, one gets a short exact sequence
$0\to A(3)\to \C_3(A_{2})\to \overline{B}^3\to 0$
in $\A$.
Since both $A(3)$ and $\overline{B}^3$ are in ${\sf{X}}$,
we conclude that $\C_3(A_{2})$ is in ${\sf{X}}$ as well.

Define $h_{2}: A_2 \to X$ and $k_{2}: B_2 \to A_2$ as follows:
\begin{equation*}
 \quad \, h_{2}(i)=
  \begin{cases}
    h(i)               & \text{if $i\in\{1,2,4\}$}\,\\
   \\
    (h(3),0)    & \text{if $i=3$}
    \end{cases}
 \ \ \    k_{2}(i)=
  \begin{cases}
    k(i)               & \text{if $i\in\{1,2,4\}$}\,\\

   \left( \begin{array}{cc}
                                     k(3) & 0 \\
                                                      0 & 1 \\
                                                    \end{array}
                                                  \right) & \text{if $i=3$}
  \end{cases}
\end{equation*}

According to the above constructions,
one gets an exact sequence
$$\bE_{2}:\,\,\,0\to B_{2} \overset{k_{2}}\longrightarrow A_{2}\overset{h_{2}}\longrightarrow X\to0$$
in $\Rep{Q}{\A}$, which is in the following form:

$$\xymatrix@R=0.005cm{
  &&& 0 \ar[ddd] \\
  &0\ar[ddd]&&&& 0 \ar[ddd] &  0 \ar[ddd]\\
  &&&& 0 \ar[ddd]\\
  &&& B(1)\ar[drr]^{\quad \quad B_2(a)}\ar[ddd]_{k(1)} \\
  &~B_2:\ar[ddd]_{k_{2}}&&&& B(3)\oplus D_2^3 \ar[r]^{\quad B_2(c)} \ar[ddd]_{k_2(3)} &  B_2(4)\ar[ddd]_{k_2(4)} \\
  &&&& B(2)  \ar[ur]_{B_2(b)} \ar[ddd]_{ k(2)}\\
  &&& A(1) \ar[drr]^{\quad \quad A_2(a)} \ar[ddd]_{h(1)}\\
  &~A_2:\ar[ddd]_{h_2}&&&& A(3)\oplus D_2^3 \ar[r]^{\quad A_2(c)} \ar[ddd]_{h_2(3)} &  A(4)\ar[ddd]_{h_2(4)} \\
  &&&& A(2)  \ar[ur]_{A_2(b)\quad}\ar[ddd]_{h(2)}\\
  &&& X(1) \ar[drr]^{\quad \quad X(a)} \ar[ddd]\\
  &~X:\ar[ddd]&&&& X(3) \ar[r]^{\quad X(c)} \ar[ddd] &  B(4)\ar[ddd]\\
  &&&& X(2)  \ar[ur]_{X(b)}\ar[ddd]\\
  &&&0  \\
  &0&&&& 0  &  0 \\
  &&&&0}$$

Now, all $\varphi^{A_2}_1$, $\varphi^{A_2}_2$ and $\varphi^{A_2}_3$ are monomorphisms
with $\C_1(A_2)$, $\C_2(A_2)$ and $\C_3(A_2)$ in $\sf{X}$.
But $\varphi^{A_2}_4$ may fail to be a monomorphism
and $\C_4(A_2)$ may fail to be in ${\sf{X}}$.
However, one can process similarly to repair $\bE_{2}$ at the vertex 4; we leave it to the reader.\qed
\end{bfhpg*}


\bibliographystyle{amsplain-nodash}

\def\cprime{$'$}
  \providecommand{\arxiv}[2][AC]{\mbox{\href{http://arxiv.org/abs/#2}{\sf
  arXiv:#2 [math.#1]}}}
  \providecommand{\oldarxiv}[2][AC]{\mbox{\href{http://arxiv.org/abs/math/#2}{\sf
  arXiv:math/#2
  [math.#1]}}}\providecommand{\MR}[1]{\mbox{\href{http://www.ams.org/mathscinet-getitem?mr=#1}{#1}}}
  \renewcommand{\MR}[1]{\mbox{\href{http://www.ams.org/mathscinet-getitem?mr=#1}{#1}}}
\providecommand{\bysame}{\leavevmode\hbox to3em{\hrulefill}\thinspace}
\providecommand{\MR}{\relax\ifhmode\unskip\space\fi MR }
\providecommand{\MRhref}[2]{%
  \href{http://www.ams.org/mathscinet-getitem?mr=#1}{#2}
}
\providecommand{\href}[2]{#2}

\end{document}